\documentclass[twoside,leqno]{article}

\usepackage[letterpaper]{geometry}

\usepackage{siamproceedings}

\usepackage[T1]{fontenc}
\usepackage{amsfonts}
\usepackage{graphicx}
\usepackage{enumitem}
\usepackage{tikz}
\usepackage{float}
\usetikzlibrary{arrows.meta,shapes}

\definecolor{darkred}{RGB}{139,35,35}
\definecolor{darkgreen}{RGB}{0,128,0}
\definecolor{darkorange}{RGB}{220,150,0}
\definecolor{violet}{RGB}{128,0,128}

\newsiamremark{remark}{Remark}

\newcommand{\R}{\mathbb{R}}
\newcommand{\T}{\mathbb{T}}

\newcommand{\C}{\mathbb{C}}
\newcommand{\E}{\mathbb{E}}
\newcommand{\Prob}{\mathbb{P}}
\newcommand{\bZ}{\mathbf{Z}}
\newcommand{\cD}{\mathcal{D}}
\newcommand{\cH}{\mathcal{H}}
\newcommand{\eps}{\varepsilon}

\title{\Large On particles and waves}
\author{Carlos E.\ Kenig\thanks{Department of Mathematics, University of Chicago, Chicago, IL 60637, USA.}}
\date{}

\begin{document}

\raggedbottom

\maketitle

\begin{abstract}
We describe the works of Yu Deng, Zaher Hani and Xiao Ma on the long-time derivation of the fluid equations from the statistics of Newtonian particle collisions and briefly the earlier work of Yu Deng and Zaher Hani on the long-time derivation of the wave kinetic equation from the statistics of wave interactions in the random cubic nonlinear Schr\"odinger equation.
\end{abstract}

\textbf{Keywords.} Particles, waves, statistics, Boltzmann equation, fluid equations, wave kinetic equations, long-time derivation

\textbf{MSC Codes.} 35Q20, 35Q35, 35Q55, 35Q70, 35Q82

% Copyright Statement
\fancyfoot[R]{\scriptsize{Copyright \textcopyright\ 20XX by SIAM\\
Unauthorized reproduction of this article is prohibited}}

%% ============================================================
%%  PREFACE
%% ============================================================

\section*{Preface}

\begin{center}
\textbf{\Large Congratulations to Yu Deng!}
\end{center}

\bigskip

This note is a tribute to Yu Deng and his deep mathematics, for which he was awarded a Fields Medal at ICM 2026. Yu Deng's work is characterized by its profound application of ``bare-hands'' techniques to tackle successfully (with his collaborators) some of the most difficult and long-standing open problems at the intersection of mathematical physics, partial differential equations and probability theory. I have concentrated here on his works (with Zaher Hani and Xiao Ma, \cite{DHM1}, \cite{DHM2}) on the derivation of the Boltzmann equation for the case of a rarefied hard sphere gas, from Newtonian particles, which solves a 120 years well known conjecture, which goes back to Hilbert's 6\textsuperscript{th} problem. I also touch upon Yu Deng's works (with Zaher Hani) on the wave kinetic theory (\cite{DH1}, \cite{DH2}, \cite{DH3}). These works precede the works on particles. The last one of them \cite{DH3}, derives the wave kinetic equation in the large, from the cubic nonlinear Schr\"odinger equation, a spectacular feat also recognized by the Fields Medal. Besides its own very high intrinsic value, this work also proved to be crucial for the later work on particles. The third area in the Fields Medal citation of Yu Deng is the probabilistic well-posedness theory for nonlinear Schr\"odinger and wave equations. I have not touched upon this work in this note, so as not to exhaust the readers and myself. In this general area, Yu Deng's collaborators were Bj\"orn Bringmann, Andrea Nahmod and Haitian Yue. Their work developed the probabilistic well-posedness theory for nonlinear dispersive equations in higher dimensions, extending Bourgain's pioneering work in 1 and 2 dimensions. Yu Deng and his collaborators were able to bring the probabilistic dispersive theory to the same level of regularity as the corresponding work of Hairer and Gubinelli in the parabolic setting. This is another landmark accomplishment. See \cite{Exp} for an expository account of these works.

In writing these notes I have borrowed expository material from the papers \cite{DHM1}, \cite{DHM2} and \cite{DH3}, with permission from the authors. I am very grateful for this help. I have tried to make this very deep material as accessible as possible. I apologize for any failures in carrying this out.

\bigskip

\begin{flushright}
Carlos E.\ Kenig\\
August 2026
\end{flushright}

%% ============================================================
%%  SECTION 1: INTRODUCTION
%% ============================================================

\section{Introduction}

Mathematics and physics share a symbiotic relationship: mathematics provides the rigorous language, formulas, and logical framework that physicists use to model the universe, while physics continually drives mathematical innovation by posing new problems that require entirely new tools to solve.

Mathematical analysis developed in the 17\textsuperscript{th} century, but many of its ideas can be traced back earlier. Its modern foundations began when Fermat and Descartes developed analytic geometry, the precursor of modern calculus. A few decades later, Newton and Leibniz developed calculus independently. Newton invented calculus to describe the laws of motion and gravity. The field of analysis grew with the stimulus of many applications of calculus to the physical sciences, through differential equations, and to many other areas of science and technology. On the occasion of the International Congress of Mathematicians (ICM) in 1900, held in Paris, D.\ Hilbert presented a collection of 10 unsolved mathematical problems, with the intention of guiding mathematical research in the new century. The full list, consisting of 23 problems, was published later, first in German in the \emph{Archiv der Mathematik und Physik}, and then, in 1902, an authorized English translation (by Mary Frances Winston Newson) was published in the \emph{Bulletin of the American Mathematical Society} 8 (10): 437--479. Problem 6 in Hilbert's Paris presentation was ``Mathematical treatment of the Axioms of Physics''. This is, of course, a very broad question, which could encompass most of modern mathematical physics. In the version that Hilbert published in 1902, he narrowed down the scope of the question to: ``The investigations on the foundations of geometry suggest the problem\ldots'' ``\ldots to treat in the same manner, by means of axioms, those physical sciences in which mathematics plays an important part, in the first rank are the theory of probability and mechanics.''

The axiomatization of the theory of probability was carried out successfully by Kolmogorov in the 1930's. On the other hand, in Hilbert's 1902 paper mentioned earlier, he wrote: ``Important investigations by physicists on the foundations of mechanics are at hand; I refer to the writings of Mach, Hertz, Boltzmann and Volkmann. It is therefore very desirable that the discussion of the foundations of mechanics be taken up by mathematicians also. Thus, Boltzmann's work on the principles of mechanics suggests\ldots the problem of developing mathematically the limiting processes, there merely indicated, which lead from the atomistic view to laws of motion of continua''. This is now frequently and informally referred to as Hilbert's 6\textsuperscript{th} problem.

The breakthrough works of Deng--Hani--Ma \cite{DHM1}, \cite{DHM2} address successfully the limiting processes discussed by Hilbert, more than 120 years ago, for the case of a rarefied hard sphere gas. Our aim in this note is to describe the lengthy, involved and deep procedure developed in \cite{DHM1}, \cite{DHM2}, and to briefly discuss the earlier work \cite{DH3} of Deng--Hani, which proved to be vital for the works \cite{DHM1}, \cite{DHM2}. See Section~3 for this.

Boltzmann introduced his kinetic equation in 1872. This equation sought to give an effective macroscopic description of the statistics of a large system of interacting particles. A puzzling feature of Boltzmann's theory, which was quite controversial, is the fact that the mesoscopic description is given by a time irreversible equation (now known as the Boltzmann equation) with an entropy functional that increases forward in time, which is claimed to be derived from a microscopic system which satisfies Newton's laws of motion, which are time reversible. This is known as the ``emergence of the arrow of time''.

A more precise formulation of the question posed by Hilbert in his 1902 paper, is to give a mathematically rigorous derivation for the macroscopic equations of fluids, such as the Euler and Navier--Stokes equations, from the Newtonian laws governing microscopic particle interactions. What is now known as the Hilbert--Grad approach to this involved the derivation of Boltzmann's kinetic equation as a crucial intermediate step. This requires justifying two limits:
\begin{enumerate}[label=\roman*)]
\item The kinetic limit, where one passes from the (microscopic) Newtonian dynamics for an $N$-particle system to the Boltzmann kinetic equation for the one-particle correlation function (which reflects the statistics of the $N$-particle system, and will be defined below) in the limit as $N \to \infty$, and
\item The hydrodynamic limit, where one passes from Boltzmann's kinetic equation to the (macroscopic) equations of fluid motion (the Euler and Navier--Stokes equations) when the collision rate goes to infinity.
\end{enumerate}

The works \cite{DHM1} and \cite{DHM2} give a rigorous proof of both limiting procedures for the case of a rarefied hard sphere gas.

There is another approach to Hilbert's problem, known as Morrey's program, which aims to connect the microscopic and macroscopic sides directly. This is mathematically more demanding than the Hilbert--Grad approach carried out in \cite{DHM1}, \cite{DHM2}. Progress in the Morrey program has been made so far under additional assumptions, such as ``artificial random forcing''. See, for instance, the works of Kac \cite{Kac}, Olla--Varadhan--Yau \cite{OVY}, Quastel--Yau \cite{QY}, Mischler--Yau \cite{MY}, etc.

Returning to the Hilbert--Grad program, as carried out in \cite{DHM1}, \cite{DHM2}, we will only discuss the hydrodynamic limit briefly here. This has been studied in earlier works (\cite{BU}, \cite{Caflish}, \cite{DEL}, \cite{GSRT}, \cite{Nishida}, etc). We will restrict ourselves to the physically relevant case of three dimensions, focusing on the derivation of Boltzmann's equation, starting from Newton's laws, which are taken as axioms. This was first obtained in \cite{DHM1}, for the case of $\R^3$. It was then extended to the case of a ``box'', the torus $\T^3$, in \cite{DHM2}, where the authors also carried out the hydrodynamic limit, using the earlier results in \cite{BU}, \cite{Caflish}, \cite{DEL}, \cite{GSRT}, \cite{Nishida}, etc.

Thus, when considering a microscopic system formed of $N$ particles of diameter $\eps$, undergoing elastic collisions, in \cite{DHM1}, \cite{DHM2} (for the case of $\R^3$, $\T^3$ respectively), the authors take the kinetic limit as $N \to \infty$, $\eps \to 0$. They show that the one-particle correlation function of the particle system (which reflects the statistics of the particle system) is a good approximation to the solution $f(t,x,v)$ of the Boltzmann equation
\begin{equation}\label{eq:boltzmann_intro}
(\partial_t + v \cdot \nabla_x) f(t,x,v) = \alpha\, \mathcal{Q}(f), \qquad f(0) = f_{\mathrm{in}},
\end{equation}
where $\alpha$ is the collision rate and $\mathcal{Q}(f)$ is the hard sphere collision kernel defined by
\begin{align*}
\mathcal{Q}(f)(t,x,v) &= \int_{(\R^3)^3} \delta(v - v_1 + v_2 - v_3)\, \delta(|v|^2 - |v_1|^2 + |v_2|^2 - |v_3|^2) \\
&\qquad\qquad \cdot (f_1 f_3 - f f_2)\, dv_1\, dv_2\, dv_3,
\end{align*}
where $f_j = f(t,x,v_j)$, $v_1, v_2, v_3 \in \R^3$, $t > 0$, $x \in \R^3$, and $\delta$ is the Dirac delta.

Here $\alpha = N\eps^2$ stands for the collision rate of the particle system, which is kept constant in the kinetic limit. The necessity of this scaling relationship between $N$ and $\eps$ (which corresponds to diluted gas) was discovered by Grad \cite{Grad}; it is referred to as the Boltzmann--Grad limit. In the hydrodynamic limit, one derives the equations of fluid dynamics when the collision rate $\alpha$ tends to infinity. Establishing the link between the two limits requires deriving the Boltzmann equation on time intervals of length $O(1)$, which corresponds (rescaling $x$ and $t$) to solutions of the (collision rate equal to 1) Boltzmann equation that exist on time intervals of length at least $\alpha$, where $\alpha$ is the original collision rate. Since $\alpha \to \infty$ in the hydrodynamic limit, it is then necessary to obtain a long-time derivation of the Boltzmann equation, which requires the main results in \cite{DHM1}, \cite{DHM2}. Obtaining the kinetic limit has turned out to be significantly more challenging than obtaining the hydrodynamic limit. In fact, the formulation of the kinetic limit as a clear mathematical question was not possible until Grad \cite{Grad} specified the precise Boltzmann--Grad scaling law, $N\eps^2 = O(1)$, in order for the Boltzmann equation to be derived from Newton's laws as $N \to \infty$, $\eps \to 0$.

Historically, after the pioneering works of Grad \cite{Grad}, and Cercignani \cite{Cerc}, the first major development in the derivation of the Boltzmann equation from Newton's laws was due to Lanford \cite{Lanford} whose breakthrough first completed the derivation for short time. Boltzmann believed that his 1872 equation followed from Newton's laws, but only if particles don't collide repeatedly, thus preserving the initial independence. This was at the core of Lanford's work, who showed this for short times (see also \cite{King}). However for larger times there may be many re-collisions. This meant that to make further progress, it was necessary to analyze a large number of collision patterns, which threaten the independence of the particles. This was a main contribution of \cite{DHM1} and \cite{DHM2}. After Lanford's breakthrough, many works were produced (\cite{BGSR1}, \cite{BGSR2}, \cite{BGSR3}, \cite{GSR2}, \cite{BGSR4}, \cite{BGSR5}, \cite{GSR}, \cite{GSRT}, \cite{IP1}, \cite{IP2}, \cite{PS}). Nevertheless all of these results were restricted to short times, or small solutions or linearized problems. These restrictions prevented the full realization of Hilbert's program until the works \cite{DHM1}, \cite{DHM2}. In \cite{DHM1}, where the setting is $\R^3$, the authors noticed that there is a ``simplification'', because particles have room to disperse and thus eventually stop colliding at all. After that, in \cite{DHM2}, where they were able to obtain the corresponding result in a finite box, particles are increasingly likely to re-collide, but they could show that the impact of multiple collisions on the independence of the system was statistically negligible.

%% ============================================================
%%  SECTION 2: PARTICLES (skeleton — to be filled in)
%% ============================================================

\section{Particles}

\subsection{Newton's Laws}

We start by a description of Newton's laws governing microscopic particle interactions (hard sphere dynamics) under the Boltzmann--Grad scaling law. This is taken directly from 1.1.1.\ in \cite{DHM1}. We thus consider $N$ particles in 3 dimensional space $\R^3$, with each particle being a ``hard sphere'' of diameter $\eps$. Here $N \to \infty$, $\eps \to 0$, and we assume the Boltzmann--Grad scaling relation $N\eps^2 = \alpha$, where $\alpha$ is the ``average collision rate''. This corresponds to the physical scenario of rarefied gas. The description of the hard sphere dynamics is as follows: let $(x_i, v_i) = z_i$, $\bZ_N = (z_i)_{i=1}^N$, $x_i \in \R^3$, $v_i \in \R^3$ represent position and velocity. $z_j$ and $\bZ_N$ are called the state vector of the $j$'th particle and the collection of all $N$ particles respectively.

The following equations hold:
\begin{equation}\label{eq:Newton}
\begin{cases}
\partial_t x_i = v_i, \quad \partial_t v_i = 0 & \text{if } |x_i - x_j| > \eps, \\[6pt]
v_i^{\mathrm{post}} = v_i^{\mathrm{pre}} - [(v_i^{\mathrm{pre}} - v_j^{\mathrm{pre}}) \cdot \omega]\, \omega, \\
v_j^{\mathrm{post}} = v_j^{\mathrm{pre}} - [(v_i^{\mathrm{pre}} - v_j^{\mathrm{pre}}) \cdot \omega]\, \omega, \\
\quad \text{if } |x_i - x_j| = \eps, \text{ where } \\
\quad \omega = \eps^{-1}(x_i - x_j) \text{ when } |x_i - x_j| = \eps,
\end{cases}
\end{equation}
$v_i^{\mathrm{post}}$ and $v_i^{\mathrm{pre}}$ denote limits from the right and left as $t$ tends to the collision time.

After some technical preliminaries this can be defined as a time reversible and volume preserving Hamiltonian dynamical system. We define the non-overlapping domain
\[
\cD_N = \bigl\{ \bZ_N = (z_1, \ldots, z_N) \in \R^{6N} : |x_i - x_j| > \eps \quad \forall\, i \neq j \bigr\}.
\]
We have $\bZ_N(0) = \bZ_N^0$. The $x_j(t)$ are always continuous in $t$; we also require $v_j(t)$ to be left continuous in $t$, so that $v_j(t) = v_j(t)^-$. If for a certain $i$ there is no $j$ such that the scenario of \cref{eq:Newton} holds, we then have $\frac{d}{dt}(x_i, v_i) = (v_i, 0)$. We define the flow map $\cH_N(t)$ by
\[
\cH_N(t) : \cD_N \longrightarrow \cD_N, \qquad \cH_N(t)(\bZ_N^0) = \bZ_N(t),
\]
where $\bZ_N(t)$ is defined by \cref{eq:Newton}. For a function $\mathfrak{f} = \mathfrak{f}(\bZ_N)$, we define the flow operator $S_N(t)(\mathfrak{f})(\bZ_N) = \mathfrak{f}(\cH_N(t)^{-1} \bZ_N)$. We define a collision when \cref{eq:Newton} holds. These collisions correspond exactly to discontinuities of $v_j(t)$. (See \cite{DHM1}, 1.1.1.\ for further details).

\textbf{\underline{Basic properties}:} Up to a set $\mathcal{Z}$ of measure 0 in $\cD_N$, we have that the hard sphere system above exists, is unique, and satisfies:
\begin{enumerate}[label=(\alph*)]
\item No two collisions happen at the same time.
\item The total number of collisions has an upper bound that depends only on $N$.
\item The flow maps $\cH_N(t)$ are measure preserving diffeomorphisms from $\cD_N \setminus \mathcal{Z}$ to itself and satisfy the semi-group property
\[
\cH_N(t+s) = \cH_N(t)\, \cH_N(s), \qquad t, s \geq 0.
\]
\end{enumerate}

\begin{figure}[H]
\centering
\begin{tikzpicture}[>=Stealth, thick]
\tikzset{
  circ/.style={circle, minimum size=0.3cm, draw=black, fill=white},
}

% Coordinates — wider gap so bar is visible between circles
\coordinate (center_i) at (-0.25, 0);
\coordinate (center_j) at (0.25, 0);
\coordinate (vi_in) at (-2.0, -2.8);
\coordinate (vj_in) at (2.0, -2.8);
\coordinate (vi_out) at (-2.0, 2.8);
\coordinate (vj_out) at (2.0, 2.8);

% Trajectory lines with arrows (colored)
\draw[darkred, line width=1.2pt, -{Stealth[length=7pt, width=5pt]}]
    (vi_in) -- (center_i);
\draw[darkgreen, line width=1.2pt, -{Stealth[length=7pt, width=5pt]}]
    (vj_in) -- (center_j);
\draw[darkred, line width=1.2pt, -{Stealth[length=7pt, width=5pt]}]
    (center_i) -- (vi_out);
\draw[darkgreen, line width=1.2pt, -{Stealth[length=7pt, width=5pt]}]
    (center_j) -- (vj_out);

% Nodes (circ style)
\node[circ] at (vi_in) {};
\node[circ] at (vj_in) {};
\node[circ] at (vi_out) {};
\node[circ] at (vj_out) {};
\node[circ] at (center_i) {};
\node[circ] at (center_j) {};

% epsilon-omega: black bar drawn LAST (on top of circles)
\draw[black, line width=2.5pt] (center_i) -- (center_j);

% Labels
\node[left, font=\normalsize] at (-1.7, -1.4) {$v_i(t)$};
\node[right, font=\normalsize] at (1.7, -1.4) {$v_j(t)$};
\node[left, font=\normalsize] at (-1.7, 1.4) {$v_i(t^+)$};
\node[right, font=\normalsize] at (1.7, 1.4) {$v_j(t^+)$};
\node[above, font=\normalsize] at (0, 0.35) {$\varepsilon\omega$};
\end{tikzpicture}
\caption{Hard sphere collision.}
\label{fig:collision}
\end{figure}
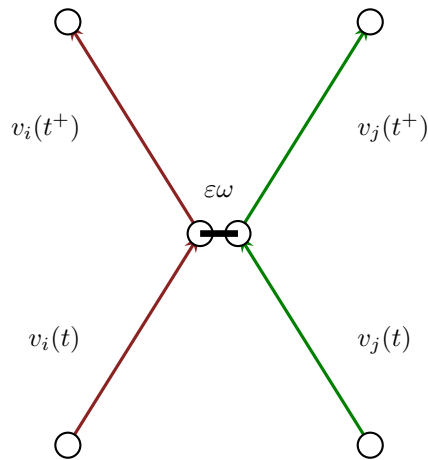

In this illustration $(v_i(t), v_j(t)) = (v_i(t^-), v_j(t^-))$ and $(v_i(t^+), v_j(t^+))$ are incoming and outgoing velocities. $\varepsilon\omega$ is the vector connecting the centers of the two colliding particles, which has length $\varepsilon$.

Since, as was predicted by Maxwell and Boltzmann in the period 1860--1870's, it is the ``one particle correlation function'' (to be defined later) of the particle system, that yields a solution of the Boltzmann equation, we must introduce statistics into the particle system. We need a formal description of that, which will be done in 2.2. This material is taken from 1.1.1.\ and 1.1.2.\ in \cite{DHM1}, where further details can be found. On first reading, the reader could just browse over 2.2 and move on.

\subsection{The grand canonical ensemble}

We define $\cD = \bigcup_N \cD_N$ to be the grand canonical domain, so $\bZ = \bZ_N \in \cD_N$ for some $N$. We define the hard sphere dynamics $\bZ(t)$ on $\cD$ and the flow map $\cH(t)$ on $\cD$ as $\bZ(t) = \bZ_N(t)$, with $\bZ^0 = \bZ_N^0$.

\textbf{\underline{Random data and initial density}}: Fix $0 < \eps \leq 1$, and $f_0(z) \geq 0$, $z = (x,v)$, $x \in \R^3$, $v \in \R^3$, $\int f_0(z)\, dz = 1$. Assume that $\bZ^0 \in \cD$ is the random variable whose law is given by the initial density functions $(W_{0,N})$, in the sense that
\[
\Prob\bigl(\bZ^0 = \bZ_N^0 \in A \subseteq \cD_N\bigr) = \frac{1}{N!} \int_A W_{0,N}(\bZ_N)\, d\bZ_N,
\]
for any $N$ and $A \subseteq \cD_N$, where $W_{0,N}$ is given by
\[
\frac{1}{N!}\, W_{0,N}(\bZ_N) = \frac{1}{\mathcal{Z}}\, \frac{\eps^{-2N}}{N!} \prod_{j=1}^N f_0(z_j)\, \mathbf{1}_{\cD_N}(\bZ_N),
\]
and where the partition function $\mathcal{Z}$ is given by
\[
\mathcal{Z} = 1 + \sum_{N=1}^{\infty} \frac{\eps^{-2N}}{N!} \int_{\R^{6N}} \prod_{j=1}^N f_0(z_j)\, \mathbf{1}_{\cD_N}(\bZ_N)\, d\bZ_N,
\]
and $\mathbf{1}_B$ is the characteristic function of the set $B$.

Let $\bZ^0 \in \cD$ be the random variable defined from $f_0$ as above, and let $\bZ(t) = \cH(t)\bZ^0$ be the evolution of the initial datum $\bZ^0$ by the hard sphere dynamics. Then, $\bZ(t)$ is also a $\cD$ valued random variable, whose density function $W_N(t, \bZ_N)$ for the law of the random variable $\bZ(t)$ is given by
\[
\Prob\bigl\{\bZ(t) = \bZ_N(t) \in A \subset \cD_N\bigr\} = \frac{1}{N!} \int_{\R^{6N}} W_N(t, \bZ_N)\, d\bZ_N.
\]
By the volume preserving property of $\cH_N(t)$, we have
\[
W_N(t, \bZ_N) = (S_N(t) W_{0,N})(\bZ_N).
\]
We can now define the $s$-particle correlation function by
\[
f_s = f_s(t, \bZ_s) = \eps^{2s} \sum_{m=0}^{\infty} \frac{1}{m!} \int_{\R^{6m}} W_{s+m}(t, \bZ_{s+m})\, dz_{s+1} \cdots dz_{s+m}.
\]
The correlation functions give a way to encode the statistical properties of the particle system. Here we have abbreviated $\bZ_s = (z_1, \ldots, z_s)$ and $\bZ_{s+m} = (z_1, \ldots, z_{s+m})$ (as for $\bZ_N$). Because of the choice of the power $(-2N)$ for $\eps$ in the formula for the partition function $\mathcal{Z}$, we have that
\begin{equation}\label{eq:BG_scaling}
\E(N) \cdot \eps^2 \sim 1, \text{ up to an error } O(\eps),
\end{equation}
where $\bZ^0 \in \cD$ is the random initial data defined earlier, $N$ is the random variable (the number of particles) determined by $\bZ^0 = \bZ_N^0 \in \cD_N$, and $\E(N)$ is its expected value. Hence, on average, we are considering the hard sphere dynamics with $\eps^{-2}$ particles, i.e.\ we are in the setting of the Boltzmann--Grad scaling law $N\eps^2 \simeq 1$, discussed earlier.

Parenthetically, it is noteworthy that the $s$-particle correlation functions evolve in time via the BBGKY hierarchy (Bogoliubov--Born--Green--Kirkwood--Yvon), of the form
\[
(\partial_t + v_1 \nabla_{x_1} + \cdots + v_s \nabla_{x_s})(f_s) = C_{s,s+1}(f_{s+1}),
\]
for a suitable collision kernel $C_{s,s+1}$.

Initially, the $s$-particle correlation functions are factorized:
\[
f_s(0) = \prod_{j=1}^s f_{\mathrm{in}}(x_j, v_j),
\]
where $f_{\mathrm{in}}$ is the initial value of the 1-particle correlation function defined above. This factorization at $t = 0$, suggests the independence of the states of different particles in the limit as $N \to \infty$. This is highly non-trivial to show mathematically (although it is believed to be physically intuitive). Proving this is one of the main results in \cite{DHM1}. This fact is usually referred to as ``propagation of molecular chaos'' in the physics literature. If propagation of chaos holds, then the time $t > 0$ statistics (which is determined by the $s$-particle correlation functions) is actually determined by the one particle correlation function $f_1(t, x, v)$.

\subsection{The Boltzmann equation and the statement of the main result for particles}

The Cauchy problem for the Boltzmann equation for hard sphere collisions, with initial datum $f(0,z) = f_0(z)$, $z = (x,v) \in \R^6$, $x \in \R^3$, $v \in \R^3$ is:
\begin{equation}\label{eq:boltzmann}
\begin{cases}
(\partial_t + v \cdot \nabla_x) f = \mathcal{Q}(f), \\[4pt]
f(0) = f_0, \\[4pt]
\displaystyle \mathcal{Q}(f) = \int_{(\R^3)^3} \delta(v - v_1 + v_2 - v_3)\, \delta(|v|^2 - |v_1|^2 + |v_2|^2 - |v_3|^2) \\
\qquad\qquad\qquad \cdot [f_1 f_3 - f f_2]\, dv_1\, dv_2\, dv_3,
\end{cases}
\end{equation}
where $f_j = f(t,x,v_j)$, and $\delta$ is the Dirac delta. The right-hand side of the first equation is referred to as the collision operator.

We turn to the main result in \cite{DHM1}. For functions $g : \R^3 \times \R^3 \to \C$, we define a norm
\[
\|g(x,v)\|_{\mathrm{Bol}^\beta} = \sum_{k \in \mathbb{Z}^3} \sup_{\substack{|x - k| \leq 1 \\ v \in \R^3}} e^{\beta |v|^2} |g(x,v)|.
\]

\begin{theorem}[Theorem 1 in \cite{DHM1}]\label{thm:main}
Fix $\beta > 0$, and a non-negative $f_0 = f_0(z)$, with $\int_{\R^6} f_0(z)\, dz = 1$. Suppose that the solution $f(t,z)$ of the Boltzmann equation \cref{eq:boltzmann} exists on the time interval $[0, t_{\mathrm{fin}}]$, and that
\[
\|e^{2\beta |v|^2} f(t, x, v)\|_{L^\infty_{x,v}} \leq A < \infty, \quad \text{for all } t \in [0, t_{\mathrm{fin}}]
\]
(with no size condition on $t_{\mathrm{fin}}$).

Suppose also that $f_0$ satisfies
\[
\|f_0\|_{\mathrm{Bol}^\beta} + \|\nabla_x f_0\|_{\mathrm{Bol}^{2\beta}} \leq B_0 < \infty.
\]

Consider the 3 dimensional hard sphere system of radius $\eps$ particles, as in 2.1, 2.2, with random initial datum $\bZ^0$, given by the grand canonical ensemble in 2.2, generated by $f_0$, under the Boltzmann--Grad scaling \cref{eq:BG_scaling}. Let $\eps$ be small, depending on $(t_{\mathrm{fin}}, \beta, A, B_0)$. Then, uniformly in $t \in [0, t_{\mathrm{fin}}]$ and in $1 \leq s \leq \log 1/\eps$, the $s$-particle correlation functions $f_s(t)$, defined in 2.2, satisfy
\[
\biggl\|f_s(t, \bZ_s) - \prod_{j=1}^s f(t, z_j)\, \mathbf{1}_\cD(\bZ_s)\biggr\|_{L^1(\R^{6s})} \leq \eps^\theta,
\]
where $\theta$ is an absolute constant. Moreover, propagation of chaos holds in the full interval $[0, t_{\mathrm{fin}}]$.
\end{theorem}

\begin{corollary}\label{cor:BG}
Under the hypotheses of \cref{thm:main}, we have
\[
\|f_1(t, x, v) - f(t, (x,v))\, \chi_\cD(x,v)\|_{L^1(\R^6)} \leq \eps^\theta, \quad \text{for } t \in [0, t_{\mathrm{fin}}],
\]
where $f_1$ is the 1-particle correlation function generated by $f_0$, and $f(t, (x,v))$ is the solution of the Boltzmann equation with initial datum $f_0$.
\end{corollary}

\cref{cor:BG} verifies the Boltzmann--Grad conjecture.

\begin{remark}
1.\ The Maxwellian (Gaussian) decay of $f(t,z)$, $z = (x,v)$, in $v$ is standard for solutions of the Boltzmann equation. The required integrability in $x$ for $f_0$ is at the level of $L^1_x$, which is the minimal requirement to define and normalize the ensemble.

2.\ For $s = 1$, it is conjectured that the minimal value of $\theta$ should be $1^-$.
The limiting (in $\eps$) value of $\E(N) \cdot \eps^2$ in \cref{eq:BG_scaling} can be replaced by any $\alpha > 0$.

3.\ The corresponding result was proved in \cite{DHM1} for $\R^d$ for any $d \geq 2$, with no change in the proofs.

4.\ The corresponding result for $\T^d$, $d = 2, 3$ was proved in \cite{DHM2}.
\end{remark}

\subsection{Rough description of the proof, with comments, part 1: Cumulants}

Lanford's proof of the short time result is based on matching two different time expansions. Expansion I is given by the expansion of the dynamics of the microscopic system satisfied by the correlation functions, through the standard Duhamel formula, and a power series expansion in $t$. (Recall that if $A$ is a $t$-independent operator, and $h$ solves $\partial_t h + Ah = g$, $h(0) = h_0$, Duhamel's formula gives $h(t) = W(t) h_0 + \int_0^t W(t-s) g(s)\, ds$, where $W(t) h_0$ solves the same problem with $g = 0$. This is a version of the method of the variation of the constants). Expansion II is given by the use of Duhamel's formula for Boltzmann's kinetic operator and another power series expansion in $t$. The main part of Lanford's proof is identifying the leading terms (in terms of powers of $t$), in the expansion I, and matching them order by order, in powers of $t$, in the expansion II, and proving that all remaining terms tend to 0 as $t$ tends to 0. For large time both expansions diverge, and the method fails. To attempt to go beyond this, and reach an arbitrarily large time $\bar{t} = t_{\mathrm{fin}}$, a natural idea is to subdivide the long interval $[0, \bar{t}]$ into $L$ subintervals (called time layers) $[(l-1)\tau, l\tau]$, $1 \leq l \leq L$, with $\tau$ sufficiently small, so that in each subinterval, the solution of the kinetic equation is given by a convergent series expansion with initial time $(l-1)\tau$. This time layering in \cite{DHM1} is the first key idea that allows to resolve the divergence problem of the expansion II, which gives the solution to the Boltzmann equation. Then, the remaining task, is to control expansion I for the microscopic system on the interval $[(l-1)\tau, l\tau]$. The obvious difficulty that arises at this point is that the data at time $(l-1)\tau$ will no longer have independence between different particles. This independence (or equivalently the complete factorization of the correlation functions $f_s$) holds at time 0, and plays a fundamental role in controlling the expansion I of the microscopic system, in the short time proof. Hence, in order to go from time $(l-1)\tau$ to $l\tau$, we need to propagate some information for the ensemble, at time $l\tau$, that describes the departure from independence. This is given by cumulants. (See 1.3.2.\ in \cite{DHM1}). The cumulants of a random variable $X$ are defined using the cumulant generating function (CGF) $K(t)$ defined by $K(t) = \log \E[e^{tX}]$. The cumulants are then defined from the power series expansion of the cumulant generating function, so that $\kappa_m = K^{(m)}(0)$, where $\kappa_m$ is the $m$'th cumulant and $K^{(m)}(0)$ is the $m$'th derivative of $K(t)$ at $t = 0$. Cumulants help to simplify the study of sums of random variables, since if $X_1, \ldots, X_k$ are independent, the cumulant of their sum equals the sum of the cumulants). The use of cumulants is crucial for the proof of the results in \cite{DHM1}. Cumulant expansions have also been used in \cite{IP1}, \cite{BGSR2}, \cite{BGSR3} and \cite{BGSR4}. They have also played a central role in the wave setting \cite{DH3}. In the present context, they are used in the following expansion of the $s$-particle correlation functions:
\begin{equation}\label{eq:cumulant_expansion}
\begin{cases}
\displaystyle f_s(t, \bZ_s) = \prod_{j \in [s]} f^{\#}(t, z_j) + \sum_{\emptyset \neq H \subset [s]} \biggl(\prod_{j \in [s] \setminus H} f^{\#}(t, z_j)\biggr) \cdot E_H(t, \bZ_H) \\[10pt]
\displaystyle \quad = \sum_{H \subset [s]} \biggl(\prod_{j \in [s] \setminus H} f^{\#}(t, z_j)\biggr) \cdot E_H(t, \bZ_H),
\end{cases}
\end{equation}
where $[s] = \{1, 2, \ldots, s\}$, $H \subset [s]$ and $\bZ_H = (z_j)_{j \in H}$.

(The actual expansion used in \cite{DHM1} has an additional error term, which is extremely small and which we will ignore here). The first term on the right hand side represents the fully factorized part of the correlation function, with $f^{\#}(t, z)$ being a 1-particle correlation function, and the remaining terms, involving the cumulants $E_H$, describing the deviation from independence.

The physical interpretation of this expansion (from 1.3.2.\ of \cite{DHM1}) is as follows. The left-hand side denotes the probability density of finding $s$ particles, say $1, 2, \ldots, s$, with states $z_1, z_2, \ldots, z_s$ at time $t$. This depends on the collision history of these particles, among themselves, and with the remaining particles. The first term on the right hand side comes from the collision histories where the $s$ particles are ``disconnected'', i.e.\ they do not interact with each other and each of them interacts with a disjoint set of the remaining particles. The second term comes from collision histories where the particles in $H$ are connected with each other via collisions. The second equality is a concise version of the first equality. From the first decomposition we see that \cref{thm:main} follows if one can show:
\begin{enumerate}
\item[(1)] $f^{\#}$ converges to the Boltzmann solution, and
\item[(2)] $\|E_H\|_{L^1}$ converges to 0 as $\eps \to 0$.
\end{enumerate}
The proof of (1) is simpler, see Proposition 6.1 in \cite{DHM1}, and Section 14 of \cite{DHM1}, so we will focus on (2).

In order to prove (2), one shows that
\begin{equation}\label{eq:cumulant_bound}
\|E_H(l\tau)\|_{L^1} \leq e^{\gamma |H|}, \qquad \forall\, \gamma > 0,
\end{equation}
where $|H|$ = the number of elements in $H$, (and $H \neq \emptyset$).
For $l = 0$, \cref{eq:cumulant_bound} is obvious, for $l = 1$, this follows from \cite{PS} in their proof of Lanford's theorem. Thus, the issue is how to treat the intervals $[(l-1)\tau, l\tau]$ to prove \cref{eq:cumulant_bound} at time $l\tau$, $l > 1$. In \cite{PS}, the authors prove that if the cumulants at time $(l-1)\tau$ satisfy the $L^\infty$ estimate $\|E_H((l-1)\tau)\|_{L^\infty} \leq \eps^{\beta |H|}$, then one can recover the estimate \cref{eq:cumulant_bound} at time $l\tau$. Thus, one starts from the stronger $L^\infty$ bound at time $(l-1)\tau$ and recovers only the weaker $L^1$ bound at time $l\tau$. However this does not iterate. In fact, as is mentioned in \cite{DHM1}, one can show that such an inductive proof of \cref{eq:cumulant_bound} is not possible, due to the time reversibility of the Newtonian dynamics and the time irreversibility of the Boltzmann equation. Thus, to prove \cref{eq:cumulant_bound}, one needs to propagate more than norms, but instead ``structural information'' on the cumulants $E_H(l\tau)$, by induction on $l$. This ``structural information'' comes from a ``partial time'' expansion, in another important new element of the proof. Let's describe how this is done in \cite{DHM1}. Start from a cumulant $E_H(l\tau)$, and assume that $f^{\#}((l-1)\tau)$ has already been approximated by the Boltzmann solution. Then express $E_H(l\tau)$ as a sum of Duhamel integrals involving $f^{\#}((l-1)\tau)$ and $E_{H}((l-1)\tau)$ (this is done through a ``cluster expansion''. Clusters are subsets of particles that are connected by collisions. For details see Section 5, \cite{DHM1}). The crucial point in the ``partial time'' expansion is that you don't expand $f^{\#}((l-1)\tau)$ any further in time, but use instead its closeness to the Boltzmann solution, and only expand the cumulant terms $E_{H}((l-1)\tau)$ into $f^{\#}((l-2)\tau)$ and $E_{H}((l-2)\tau)$, avoiding the expansion of the leading terms $f^{\#}((l-1)\tau)$ all the way to time 0. This bypasses the divergence remarked on earlier, in the Lanford type arguments.

The cumulants come with power gains (of size $\eps^{\gamma |H|}$) from \cref{eq:cumulant_bound}, which offsets the divergence in the expansion in powers of $t$. This is crucial in obtaining a time irreversible equation (Boltzmann) from time reversible dynamics (Newton). The propagation of the cumulants is what sees the emergence of irreversibility.

\subsection{Rough description of the proof, with comments, part 2: Reduction to combinatorics}
Because of the physical interpretation of the cumulants and the ``partial time'' expansion outlined in 2.4, we know that the cumulants $E_H(l\tau)$ can be written as sums of contributions corresponding to probability densities of possible collision histories on $[0, l\tau]$, involving particles in $H$. These collision histories are constructed by the ``partial time'' expansion, and the particles in $H$ are connected to each other via collisions. Next, we turn to a brief description of the most crucial and novel part of the proof of \cref{thm:main}. This is a tour de force in combinatorics. At this point in the proof, \cite{DHM1} disregards the precise positions, velocities and collision times of the particle trajectories and focuses only on the combinatorial structure. Each collision history topologically (in shape) reduces to an abstract diagram, called a molecule $\mathbb{M}$. These diagrams encode the possible patterns of collisions among the particles in $H$, and disregard the precise geometry. To illustrate this we start with a 2-dimensional depiction of the trajectories of four particles, as a function of time.

% Particle trajectories and topological reduction
\begin{figure}[H]
\centering
\includegraphics[width=0.95\textwidth]{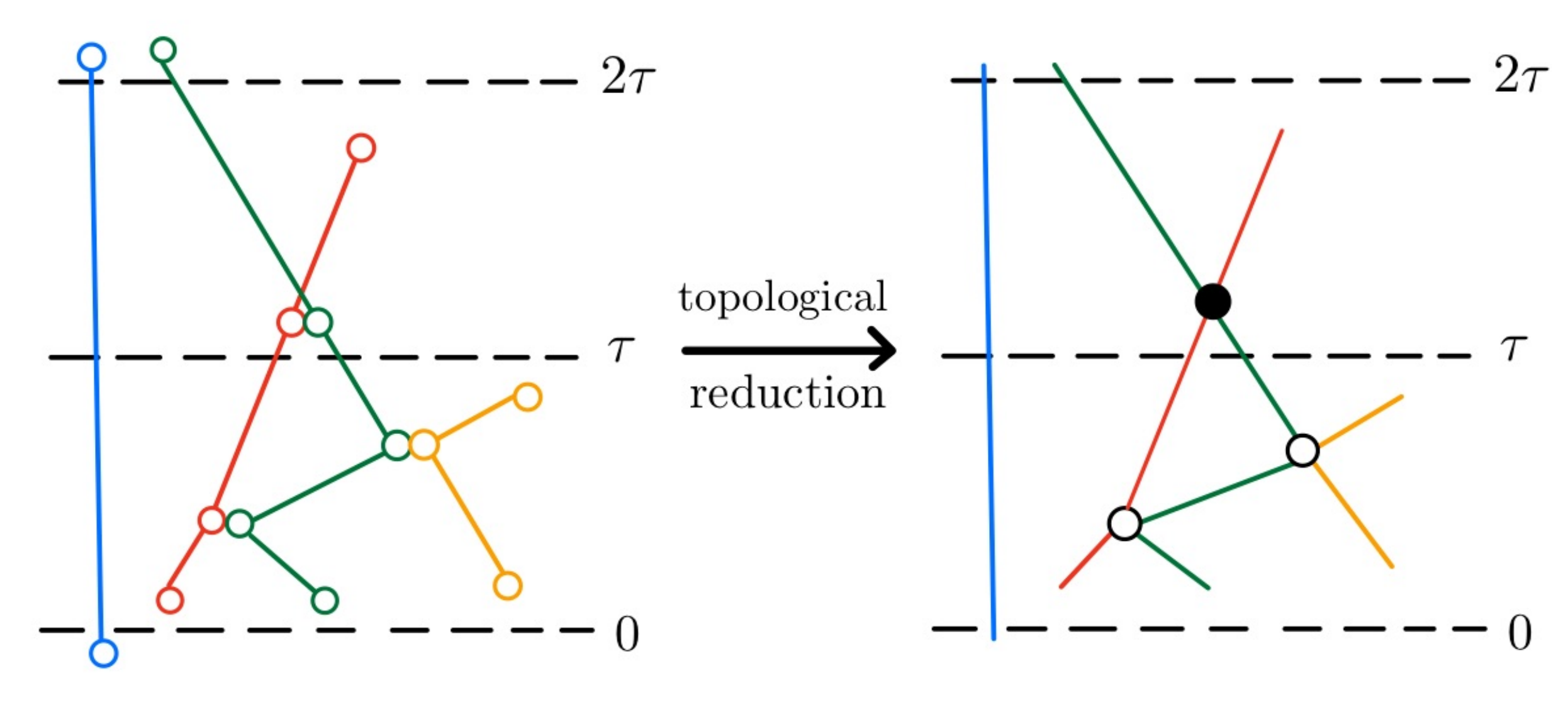}
\caption{Left: Trajectories of four particles across two time layers. Atoms (depicted by circles) denote collisions or overlaps. Right: The corresponding molecule after topological reduction --- endpoint atoms are removed, remaining overlaps are colored black (filled), and collisions are colored white (open).}
\label{fig:trajectories}
\end{figure}

The particles follow the different colored lines, atoms (depicted by circles) denote either collisions or overlaps of different particles. When we do the ``topological reduction'' we eliminate all atoms at the end of a trajectory, and we color in black all remaining overlaps and in white all collisions. This is the corresponding molecule, which appears after the ``topological reduction''.

Note also that the molecules retain the ordering between different collisions of the same particle, and the time layer of each collision, represented by the black horizontal lines.

Using (layered) molecules, the ``structural information'' for the cumulants $E_H$ can now be stated as
\begin{equation}\label{eq:molecule_bound}
|E_H(l\tau)| \leq \sum_{\mathbb{M}} |\mathcal{I}\mathcal{N}_{\mathbb{M}}|,
\end{equation}
where the sum is taken over molecules $\mathbb{M}$ constructed from the ``partial time'' expansion process developed to prove \cref{eq:cumulant_bound}, and $|\mathcal{I}\mathcal{N}_{\mathbb{M}}|$ is an explicit quantity representing the (normalized) probability density that the collisions described by the molecule $\mathbb{M}$ actually happen. To obtain the required bound on $\|E_H\|_{L^1}$, in (2) above, one has to prove effective estimates for the $L^1$ norm of the $|\mathcal{I}\mathcal{N}_{\mathbb{M}}|$ just introduced. The central quantity for understanding such estimates is the number of re-collisions $\rho$, in the molecule $\mathbb{M}$. (A re-collision is a new collision between particles that had already collided). When a re-collision takes place in $\mathbb{M}$, a ``cycle'' is created in the molecule $\mathbb{M}$. Heuristically, as the number $\rho$ of re-collisions increases, the number of molecules in the sum \cref{eq:molecule_bound} also increases. To obtain the desired estimate \cref{eq:cumulant_bound} one then has to show a gain in the estimates for $|\mathcal{I}\mathcal{N}_{\mathbb{M}}|$ for molecules with larger $\rho$. Given $\rho \geq 1$, the following estimates are proven in \cite{DHM1}:
\begin{equation}\label{eq:molecule_count}
\begin{gathered}
\text{The number of molecules } \mathbb{M} \text{ with } \rho \text{ re-collisions verifies} \\
\#\,\mathbb{M} \leq C^{|\mathbb{M}|} \cdot (\log \eps)^{C\rho}, \quad \text{for suitable } C,
\end{gathered}
\end{equation}
where $|\mathbb{M}|$ is the number of collision atoms in the molecule $\mathbb{M}$. (Overlap atoms don't affect the dynamics and can thus be ignored). The contribution of each molecule $\mathbb{M}$ is
\begin{equation}\label{eq:molecule_contribution}
\|\mathcal{I}\mathcal{N}_{\mathbb{M}}\|_{L^1} \leq \tau^{|\mathbb{M}|} \cdot \eps^{\gamma \rho}, \qquad \gamma > 0.
\end{equation}

These estimates explain the effect of re-collisions: Each additional re-collision leads to a $|\log \eps|^C$ increase in the number of molecules. This is compensated by the $\eps^\gamma$ gain in the normalized probabilities in \cref{eq:molecule_contribution}. The desired cumulant estimate follows from these properties, provided that $\eps \ll \tau \ll 1$.

The constructed molecules $\mathbb{M}$ consist of $l$-layers, and they can be seen as a concatenation of $l$ submolecules $\mathbb{M}_{l'}$, with $1 \leq l' \leq l$, and with $\mathbb{M}_{l'}$ containing the atoms in the layer $l'$. To establish the desired estimates \cref{eq:molecule_count} and \cref{eq:molecule_contribution}, one imposes an upper bound on the number of atoms and recollisions in each submolecule $\mathbb{M}_{l'}$. This leads to the notion of ``truncated dynamics'' within each time layer. (See 4.3 in \cite{DHM1} for the precise definition). In the ``truncated dynamics'', at each time where a collision would occur, one checks whether or not this collision would cause a violation to the following two conditions (recalling that a cluster is a subset of particles that are connected by collisions):
\begin{enumerate}[label=\alph*)]
\item The number of particles in each cluster is at most $\Lambda = |\log \eps|^{O(1)}$.
\item The number of re-collisions in each cluster is at most $\Gamma = O(1)$.
\end{enumerate}
If a collision would cause a violation, one ``turns off'' this collision and allows the particles to cross each other without colliding, otherwise the collision is allowed as in the original dynamics. This truncation leads to an error, but it can be controlled again in terms of $|\mathcal{I}\mathcal{N}_{\mathbb{M}}|$, and treated in a similar way as $E_H$.

For the truncated dynamics, a) easily follows from a standard fact in combinatorics, which states that the number of graphs whose vertices have maximum degree 4 (the degree of a vertex is the number of connections that the vertex has to other vertices), with at most $m$ vertices and $\rho$ independent cycles, is bounded by $C^m \cdot m^\rho$, for some universal $C$. This fact follows from the counting upper bound for binary and ternary trees. This leaves b) to be proved, which is the heart of the proof. To carry this out, \cite{DHM1} is led to constructing a ``cutting algorithm'', a purely combinatorial problem for a molecule $\mathbb{M}$. One first re-writes $\|\mathcal{I}\mathcal{N}_{\mathbb{M}}\|_{L^1} = \eps^{2|H|} \cdot \mathcal{I}_{\mathbb{M}}(Q)$, where $\mathcal{I}_{\mathbb{M}}$ is an operator involving integration over all intermediate positions $x_i$ and velocities $v_i$ associated with the edges of $\mathbb{M}$, and $Q$ is a specific non-negative function. The integrals involve $\delta$ functions, which encode the relations between incoming and outgoing velocities at each collision. To find the optimal bound corresponding to $\mathcal{I}_{\mathbb{M}}(Q)$, the key step is to find the correct order in which the variables should be integrated. This requires the introduction of the key notion of ``cutting'', which provides a systematic way of deciding the order of integration in $\mathcal{I}_{\mathbb{M}}$, to apply Fubini's theorem. Each ``cutting'' divides the molecule $\mathbb{M}$ into two separate molecules $\mathbb{M}_1$ and $\mathbb{M}_2$. By first integrating in the variables associated with $\mathbb{M}_2$, keeping the variables associated with $\mathbb{M}_1$ fixed, we have $\mathcal{I}_{\mathbb{M}} = \mathcal{I}_{\mathbb{M}_1} \circ \mathcal{I}_{\mathbb{M}_2}$. The strategy is to devise a ``cutting algorithm'' that cuts $\mathbb{M}$ into a number of small molecules $\mathbb{M}_j$, with one or two atoms. These $\mathbb{M}_j$ are called ``elementary molecules''. For them, $\mathcal{I}_{\mathbb{M}_j}$ can be calculated directly, and estimated (say in the $L^\infty$ to $L^\infty$ norm). Using the factorization one gets an upper bound for $\mathcal{I}_{\mathbb{M}}$. In this process, the most important ``elementary molecules'' are those which represent collisions or the ones that represent collisions followed by re-collisions. The first ones are called ``normal'', since integrating does not gain any power. The second ones are called ``good'', since integration gains a power of $\eps^\gamma$, for some $\gamma > 0$. There is also a third kind of ``elementary molecule'', which is called ``bad'' because they waste dimensions of integration. The main ingredient in the proof of b) is to construct a ``cutting algorithm'', that produces enough ``good'' elementary molecules. This is a main major novelty of the proof and a decisive technical ingredient introduced to prove \cref{thm:main}. The ``cutting algorithm'' designed in \cite{DHM1} has two big steps, (A) and (B). (A) reduces the multi-layer case to the two-layer case. (B) treats the two-layer case. In Section 11, 11.2, of \cite{DHM1}, the authors introduce an explicit combinatorial ``toy model'', for which they construct the ``cutting algorithm'', in a simplified two-layer, purely combinatorial setting. This is then extended to the general case. The reader is invited to go through 11.2 in \cite{DHM1} to appreciate the combinatorial heart of this amazing proof.

I hope that this short description (lifted from Sections 1 and 2 of \cite{DHM1}, with permission from the authors) gives an idea of the gigantic, bare-hands accomplishment by the authors of \cite{DHM1}, to carry out the proof of \cref{thm:main}. This will remain, through the times, as an amazing milestone in mathematical physics, powered by analysis and combinatorial arguments.

%% ============================================================
%%  SECTION 3: THE WAVE KINETIC THEORY
%% ============================================================

\section{The wave kinetic theory: a short discussion}

At the time of Hilbert's 1900 lecture, quantum physics was at its very beginning. With its advances in the 20\textsuperscript{th} century, especially the wave-particle duality, physicists began to develop a statistical theory of waves, known as the wave kinetic theory or wave turbulence. This led to the formulation of questions that correspond to Hilbert's 6\textsuperscript{th} problem, in the wave setting. The wave kinetic theory was initiated in works of Nordheim (1928), Peierls (1929) and Uehling--Uhlenbeck (1933). In these settings, particles are replaced by wave modes, collisions are replaced by nonlinear wave interactions (for example as in nonlinear Schr\"odinger equations) and the Boltzmann equation is replaced by the wave kinetic equation, to be introduced below. The wave kinetic theory is a huge subject with many scientific applications in such areas as plasma physics, oceanography, crystal thermodynamics, and others.

Mathematical progress on the wave kinetic theory has been much slower than in the Newtonian particles/Boltzmann theory, but there has been much progress in recent years.

We now briefly describe the wave kinetic setup in 3 dimensions, by comparison to the particle setup. Particles are replaced by Fourier modes in a large periodic box of size $L$. Collisions are replaced by nonlinear dispersive equations, of nonlinear strength $\alpha$. (In \cite{DH3}, this dispersive equation is the cubic nonlinear Schr\"odinger equation $(i\partial_t u + \Delta) u = \alpha |u|^2 u$, and the strength of the nonlinearity, $\alpha$, tends to 0 as $L \to \infty$). In \cite{DH3}, it is assumed that $\alpha L^\delta = 1$, for some scaling law parameter $\delta \in (0, 1]$. This is the full range of the admissible scaling laws, as was shown in \cite{DH2}. In analogy to the particle case, the Fourier modes of the initial data are given by
\[
u(0) = u_{\mathrm{in}} = L^{-3/2} \sum_k \hat{u}_{\mathrm{in}}(k)\, e^{2\pi i k \cdot x},
\]
where
\[
\hat{u}_{\mathrm{in}}(k) = \sqrt{n_{\mathrm{in}}(k)} \cdot \eta_k, \qquad k \in \mathbb{Z}_L^3 = L^{-1}\mathbb{Z}^3,
\]
$n_{\mathrm{in}}$ is a Schwartz class function (rapidly decreasing together with all of its derivatives) and $\eta_k$ are i.i.d.\ (independent identically distributed) normalized Gaussians.

The role of the one-particle correlation function is played by two-point correlations (see below for the definition). We define the kinetic time $T_{\mathrm{kin}} = \alpha^{-2}$. As in the particle case, from physics' considerations, it is expected to have ``propagation of chaos'', that is the initial time independence of the Fourier modes is preserved in the limit as $L$ tends to infinity. The two point correlation $n(\tau, k)$ is given by
\[
n(\tau, k) = \E\bigl(|\hat{u}(\tau \cdot T_{\mathrm{kin}}, k)|^2\bigr),
\]
and it is expected to satisfy the wave kinetic equation, for $k \in \R^3$,
\[
\partial_\tau n = \mathcal{C}(n), \qquad n(0) = n_{\mathrm{in}},
\]
\begin{align*}
\mathcal{C}(n)(\tau, k) &= \int_{(\R^3)^3} \delta(k - k_1 + k_2 - k_3)\, \delta(|k|^2 - |k_1|^2 + |k_2|^2 - |k_3|^2) \\
&\qquad \cdot \bigl[n_1 n_3 (n + n_2) - n n_2 (n_1 + n_3)\bigr]\, dk_1\, dk_2\, dk_3,
\end{align*}
where $n_j = n(\tau, k_j)$, $\delta$ is the Dirac delta function. This is a time irreversible equation. The collision operator $\mathcal{C}(n)$ has a very similar structure to the Boltzmann collision operator $\mathcal{Q}(f)$ from \cref{eq:boltzmann}, but the collision integral is cubic, instead of quadratic as in the case of Boltzmann, but with the same domain of integration. (This is specific to the choice of cubic (NLS) as the nonlinear dispersive equation). The main result in \cite{DH3} is:

\begin{theorem}[\cite{DH3}]\label{thm:wave}
Given any interval $[0, \tau^*]$ on which the solution $f(\tau, k)$ to the wave kinetic equation exists, we have
\[
\E\bigl(|\hat{u}(\tau \cdot T_{\mathrm{kin}}, k)|^2\bigr) = f(\tau, k) + O(L^{-\nu}),
\]
$\nu > 0$, on the full interval $[0, \tau^*]$, and ``propagation of chaos'' holds.
\end{theorem}

Before this breakthrough result, the state of the art was due to recent works of the same authors (\cite{DH1}, \cite{DH2}) which reached a small multiple of $T_{\mathrm{kin}}$, and which is the analog of the particle result in \cite{Lanford}. The result of \cite{DH3} was the first long-time one in this area, a very important advance. Many of the ideas in \cite{DH3} served as the foundation for the subsequent work \cite{DHM1} on particles. For example, the time layering idea, and the use of cumulants to deal with long-time independence, and with the ``emergence of the arrow of time'' first appear in \cite{DH3}. In \cite{DH3}, the authors accounted for behavior changing interactions of waves, by constructing a mathematical ``ledger'' that could reflect every possible pattern. Particles behave very differently than waves, but in \cite{DHM1} the authors use the same approach, to construct a ``ledger'' to record every collision history. In both \cite{DH3} and \cite{DHM1} the authors constructed a ``cutting algorithm'' to break down the complex histories into manageable pieces. Despite the many important differences, it is fair to say that the work in \cite{DHM1} and \cite{DHM2} was only possible because of the earlier work in \cite{DH3}, which is also extremely valuable in its own right.

%% ============================================================
%%  REFERENCES
%% ============================================================

\end{document}